\documentclass[amssymb,11pt]{amsart}
\usepackage{amscd,amssymb,euscript,amsthm}
\theoremstyle{plain}
\newtheorem{thm}{Theorem}[section] 

\newtheorem{prop}[thm]{Proposition}
\newtheorem{lem}[thm]{Lemma}
\theoremstyle{definition}
\newtheorem{defn}[thm]{Definition}
\theoremstyle{remark}
\newtheorem{rem}[thm]{Remark}

\numberwithin{equation}{section}
\newcommand{\defect}{\operatorname{defect}}

\def\<{\left<}
\def\>{\right>}
\def\cstar{$C^*$-algebra}
\begin{document}
\title[Free Cover]{The Free Cover of a row contraction}
\author{William Arveson}
%
%
%
\address{Department of Mathematics,
University of California, Berkeley, CA 94720}
\email{arveson@math.berkeley.edu}
\begin{abstract}
We establish the existence and uniqueness of finite 
free resolutions - and their attendant Betti 
numbers - for graded commuting $d$-tuples of 
Hilbert space operators.  Our approach is 
based on the notion of {\em free cover} of a 
(perhaps noncommutative) row contraction.  
Free covers provide a flexible replacement 
for minimal dilations that is better suited 
for higher-dimensional operator theory.  

For example, every 
graded $d$-contraction that is  
finitely multi-cyclic has a unique free cover 
of finite type -  whose kernel 
is a Hilbert module 
inheriting the same properties.  This contrasts 
sharply with what can be 
achieved by way of dilation theory 
(see Remark \ref{rem1}).  
\end{abstract}
%
%
\maketitle

\section{Introduction}\label{S:Int}

The central result of this paper 
establishes the existence and uniqueness 
of finite free resolutions for 
commuting $d$-tuples of operators acting on a common 
Hilbert space (Theorem \ref{thmB}).  
Commutativity is essential for that result, 
since finite resolutions do not exist for noncommuting 
$d$-tuples.  

On the other hand, we 
base the existence of free resolutions 
on a general notion 
of {\em free cover} that is effective in 
a broader noncommutative context.  
Since free covers have applications that go beyond 
the immediate needs  
of this paper, and since we intend to 
take up such applications elsewhere, 
we present the general version below 
(Theorem \ref{thmA}).  In the following section we 
give precise statements of these two results, we 
comment on how 
one passes from the larger noncommutative 
category to the commutative one, and we 
relate these results to previous work 
that has appeared 
in the literature.  Section \ref{S:Gen} 
concerns generators 
for Hilbert modules, in which we show that the examples 
of primary interest are properly generated.  
The next two sections are devoted to proofs of the 
main results - the 
existence and uniqueness of free covers and 
of finite free resolutions.  In Section \ref{S:examp} 
we discuss examples.   

In large part, this paper was composed during a visit to Kyoto University 
that was supported by MEXT, the Japanese Ministry of Education, 
Culture, Science and Technology.  It is a pleasure 
to thank Masaki Izumi for providing warm hospitality and a stimulating 
intellectual environment during this period.

\section{Statement of results}\label{S:results}

A {\em row contraction} of dimension $d$ is a $d$-tuple 
of operators $(T_1,\dots,T_d)$ acting on a common 
Hilbert space $H$ that has norm at most $1$ when viewed 
as an operator from $H\oplus\cdots\oplus H$ to $H$.  
A $d$-{\em contraction} is a row contraction whose operators 
mutually commute, $T_jT_k=T_kT_j$, $1\leq j,k\leq d$.  
In either case, one can view $H$ as a module over the noncommutative 
polynomial algebra $\mathbb C\langle z_d,\dots,z_d\rangle$ by way of 
$$
f\cdot \xi=f(T_1,\dots,T_d)\xi,\qquad f\in \mathbb C\langle z_1,\dots,z_d\rangle, 
$$
and $H$ becomes a {\em contractive} Hilbert module in the sense that 
$$
\|z_1\cdot\xi_1+\cdots +z_d\cdot\xi_d\|^2\leq  \|\xi_1\|^2+\cdots+\|\xi_d\|^2, 
\qquad \xi_1,\dots,\xi_d\in H.  
$$
The maps of this category are linear operators $A\in \mathcal B(H_1,H_2)$ 
that are homomorphisms of the module structure and satisfy $\|A\|\leq 1$.  
It will be convenient to refer to a Hilbert space 
endowed with such a module structure 
simply as a {\em Hilbert module}.  

Associated with every Hilbert module $H$ there is an integer invariant 
called the defect, defined as follows.  Let $Z\cdot H$ 
denote the closure of the range of the coordinate operators
$$
Z\cdot H=\{z_1\xi_1+\cdots+z_d\xi_d: \xi_1,\dots,\xi_d\in H\}^-.   
$$
$Z\cdot H$ is a closed submodule of $H$, hence the quotient 
$H/(Z\cdot H)$ is a Hilbert 
module whose row contraction is $(0,\dots,0)$.  
One can identify $H/(Z\cdot H)$ more concretely 
as a subspace of $H$ in terms of the ambient operators 
$T_1,\dots,T_d$, 
$$
H/(Z\cdot H)\sim H\ominus (Z\cdot H)=\ker T_1^*\cap\cdots\cap \ker T_d^*.  
$$

\begin{defn}\label{defStat1}
A Hilbert module 
$H$ is said to be {\em properly generated} 
if $H\ominus (Z\cdot H)$ is a generator:
$$
H= \overline{\rm span}\{f\cdot\zeta: f\in \mathbb C\langle z_1,\dots, z_d\rangle
\quad \zeta\in H\ominus(Z\cdot H)\}.  
$$
In general, the quotient Hilbert space $H/(Z\cdot H)$ is called the {\em defect space} 
of $H$ and its dimension $\dim(H/(Z\cdot H))$ is called 
the {\em defect}.  

\end{defn}

The 
defect space of a finitely 
generated Hilbert module must be finite-dimensional.   
Indeed, it is not hard to see that 
the defect is dominated by the smallest 
possible number of generators.  
A fuller discussion of properly generated Hilbert modules 
of finite defect 
will be found in  
Section \ref{S:Gen}.

The free objects of this
category are defined as follows.  Let $Z$ be 
a Hilbert space of dimension $d=1,2,\dots$ and let 
$ F^2(Z)$ be the Fock space over 
$Z$, 
$$
F^2(Z)=\mathbb C\oplus Z\oplus Z^{\otimes 2}\oplus
Z^{\otimes 3}\oplus \cdots
$$
where $Z^{\otimes n}$ denotes the full tensor product 
of $n$ copies of $Z$.  One can view $F^2(Z)$ 
as the completion 
of the tensor algebra over $Z$ in a natural Hilbert space 
norm; in turn, if we fix an orthonormal basis 
$e_1,\dots,e_d$ for $Z$ then we can define an 
isomorphism of the noncommutative polynomial 
algebra $\mathbb C\langle z_1,\dots,z_d\rangle$ 
onto the tensor algebra by sending $z_k$ to $e_k$, 
$k=1,\dots,d$.  This allows us to realize the Fock space as 
a completion of $\mathbb C\langle z_1,\dots,z_d\rangle$, 
on which the multiplication operators associated 
with the coordinates $z_1,\dots, z_d$ 
act as a row contraction.  We write this Hilbert 
module as $F^2\langle z_1,\dots,z_d\rangle$; and when there is no possibility 
of confusion about the dimension or choice of basis, we 
often use the more compact $F^2$.  

One forms free Hilbert modules of higher multiplicity 
by taking direct sums of copies of $F^2$.  More explicitly, let 
$C$ be a Hilbert space of dimension 
$r=1,2,\dots,\infty$ and consider the Hilbert space 
$F^2\otimes C$.  There is 
a unique Hilbert module structure on $F^2\otimes C$ 
satisfying 
$$
f\cdot(\xi\otimes \zeta)=(f\cdot\xi)\otimes \zeta, \qquad 
f\in\mathbb C\langle z_1,\dots,z_d\rangle,\quad \xi\in F^2,
\quad \zeta\in C, 
$$
making $F^2\otimes C$ into a properly generated 
Hilbert module of 
defect $r$.  

More generally, it is apparent that 
any homomorphism of Hilbert modules 
$A: H_1\to H_2$ induces a contraction 
$$
\dot A: H_1/(Z\cdot H_1)\to H_2/(Z\cdot H_2)
$$
that maps one defect space into the other.  

\begin{defn}\label{defStat2}
Let $H$ be a Hilbert module.  
By a {\em cover} of $H$ we mean a contractive 
homomorphism 
of Hilbert modules $A: F\to H$  
that has dense range and induces a 
unitary operator $\dot A: F/(Z\cdot F)\to H/(Z\cdot H)$ 
from one defect space onto the other.  A {\em free cover} of $H$  
is a cover $A: F\to H$ in which 
$F=F^2\langle z_1,\dots,z_d\rangle\otimes C$ 
is a free Hilbert module.  
\end{defn}

\begin{rem}[Extremal Property of Covers]
In general, if one is given a contractive homomorphism with dense
range $A: F\to H$, there is no way of relating 
the image $A(F\ominus (Z\cdot F))$ to $H\ominus(Z\cdot H)$, 
even when $A$ induces a bijection $\dot A$ of one defect space 
onto the other.  
But since a cover  is a contraction that 
induces a {\em unitary} map of defect spaces,  
it follows that 
$A$ must map $F\ominus (Z\cdot F) $ 
isometrically onto $H\ominus (Z\cdot H)$
(see Lemma \ref{lemPrf1}).  This extremal 
property is critical, leading for example to the uniqueness 
assertion of Theorem \ref{thmA} below.  
\end{rem}

It is not hard to give examples of finitely generated 
Hilbert modules $H$ that are degenerate in the sense 
that $Z\cdot H=H$ (see the proof of Proposition \ref{propGen2}), 
and in such cases,  
free covers $A: F\to H$ cannot exist when $H\neq\{0\}$.  As 
we will see momentarily,  the notion of free 
cover is effective for Hilbert modules that 
are {\em properly} generated.  
We 
emphasize that in a free cover $A: F\to H$ of a finitely 
generated Hilbert module $H$ with  $F=F^2\otimes C$, 
$$
\dim C=\defect (F^2\otimes C)=\defect H<\infty, 
$$
so that for finitely generated Hilbert modules for 
which a free cover exists, {\em the free 
module associated with a free cover must be of finite defect}.  
More generally, we say that a diagram of Hilbert 
modules 
$$
F\underset{A}\longrightarrow G\underset{B}\longrightarrow H 
$$
is {\em weakly exact} at $G$ if $AF\subset \ker B$ and 
the map $A: F\to \ker B$ defines 
a cover of $\ker B$.  This implies that $AF$ is dense in 
$\ker B$, but of course it asserts somewhat more.

Any cover $A: F\to H$ of $H$ can be converted into another one 
by composing it with a unitary automorphism of $F$ 
on the right.  
Two covers $A: F_A\to H$ and $B: F_B\to H$ 
are said to be {\em equivalent} 
if there is a unitary isomorphism of Hilbert modules 
$U: F_A\to F_B$ such that 
$B=A U$.  Notice the one-sided nature 
of this relation; in particular, 
two equivalent covers of a Hilbert module 
$H$ must have identical (non-closed) ranges.  
When combined with Proposition 
\ref{propGen1} below, the following 
result gives 
an effective characterization of 
the existence of free covers.  

\begin{thm}\label{thmA}
A contractive Hilbert module $H$ over 
the noncommutative polynomial algebra 
$\mathbb C\langle z_1,\dots,z_d\rangle$ 
has a free cover if, and only if, it is 
properly generated; and in that case   
all free covers of $H$ are equivalent. 
\end{thm}

\begin{rem}[The Rigidity of Dilation Theory]\label{rem1}
Let $H$ be a pure, finitely generated, contractive Hilbert module 
over $\mathbb C\langle z_1,\dots,z_d\rangle$ (see \cite{arvSubalgIII}).  
The methods of dilation theory lead to 
the fact that, up to unitary equivalence of Hilbert modules, 
$H$ can be realized as a quotient of a free Hilbert module
$$
H= (F^2\langle z_1,\dots,z_d\rangle\otimes C)/M 
$$
where $M$ is an invariant subspace 
of $F^2\otimes C$.  In more explicit terms, there is 
a contractive homomorphism $L: F^2\otimes C\to H$ of 
Hilbert modules such that $LL^*=\mathbf 1_H$.  
When such a realization is minimal, there is 
an appropriate sense in which it is unique.  

The problem with this realization of $H$ as 
a quotient of a free Hilbert module 
is that the coefficient space $C$ is 
often infinite-dimensional; moreover, the 
connecting map $L$ is only rarely a cover.  Indeed, in order for 
$C$ to be finite-dimensional it is necessary and sufficient 
that the ``defect operator" of $H$, namely 
\begin{equation}\label{EqStat0}
\Delta=(\mathbf 1_H-T_1T_1^*-\cdots-T_dT_d^*)^{1/2}, 
\end{equation}
should be of finite 
rank.  The fact is that this 
finiteness condition often fails,  
even when the underlying operators 
of $H$ commute.  

For example, any invariant subspace 
$K\subseteq H^2$ of the rank-one free commutative 
Hilbert module $H^2$, 
that is also invariant under the gauge group $\Gamma_0$
(see the following paragraphs),  
becomes a finitely 
generated graded Hilbert module whose operators 
$T_1,\dots, T_d$ are the restrictions of the $d$-shift to $K$.  
However, the defect operator of 
such a $K$ is of infinite rank in every nontrivial 
case - namely, whenever $K$ is nonzero and of infinite codimension 
in $H^2$.  Thus, even though dilation theory provides a 
realization of $K$ as the quotient 
of another free commutative Hilbert module 
$K\cong (H^2\otimes C)/M$, the 
free Hilbert module $H^2\otimes C$ must have  
{\em infinite defect}.  
\end{rem}

One may conclude from these observations that dilation 
theory is too rigid to provide an effective representation 
of finitely generated Hilbert modules as quotients of 
free modules of finite defect, and 
a straightforward application of dilation theory cannot 
lead to finite free resolutions in multivariable operator 
theory.   Our purpose below is to initiate an approach to 
the existence of free resolutions that is based on 
free covers.

We first discuss grading in the 
general noncommutative context.  
By a {\em grading} on a Hilbert module $H$ 
we mean a strongly continuous unitary representation of 
the circle group $\Gamma: \mathbb T\to \mathcal B(H)$ that 
is related to the ambient row contraction as follows 
\begin{equation}\label{EqStat1}
\Gamma(\lambda)T_k\Gamma(\lambda)^*=\lambda T_k,
\qquad \lambda\in \mathbb T,\quad k=1,\dots,d.  
\end{equation}
Thus we are restricting ourselves to gradings in which 
the given operators $T_1,\dots, T_d$ are all of degree one.  
The group $\Gamma$ is called the {\em gauge group} 
of the Hilbert module $H$.  Note that while 
there are many gradings of $H$
that are consistent with its module structure, 
when we refer to $H$ as a graded Hilbert module it is 
implicit that a particular 
gauge group has been singled out.  A graded morphism 
$A: H_1\to H_2$ of graded Hilbert modules is a 
homomorphism $A\in\hom(H_1,H_2)$ 
that is of degree zero in the sense that 
$$
A\Gamma_1(\lambda)=\Gamma_2(\lambda)A,
\qquad \lambda\in\mathbb T, 
$$
$\Gamma_k$ denoting the gauge group of $H_k$.

The gauge group of $F^2$ 
is defined by 
$$
\Gamma_0(\lambda)=\sum_{n=0}^\infty \lambda^n E_n
$$
where $E_n$ is the projection onto $Z^{\otimes n}$.  
In this way, {\em $F^2$ becomes a graded contractive Hilbert 
module over $\mathbb C\langle z_1,\dots,z_d\rangle$ 
of defect $1$}.  More generally, let $F=F^2\otimes C$ be 
a free Hilbert module of higher defect.  Since the 
ambient operators $U_1,\dots,U_d$ of $F^2$ generate 
an irreducible \cstar, one readily verifies that 
the most general strongly continuous unitary representation 
$\Gamma$ of the circle group on $F$ that satisfies 
$\Gamma(\lambda)(U_k\otimes\mathbf 1_C)\Gamma(\lambda)^*=\lambda U_k\otimes\mathbf 1_C$ 
for $k=1,\dots, d$ must decompose into a tensor 
product of representations 
$$
\Gamma(\lambda)=\Gamma_0(\lambda)\otimes W(\lambda),
\qquad \lambda\in\mathbb T,
$$
where $W$ is an arbitrary strongly continuous unitary 
representation of $\mathbb T$ on the coefficient space 
$C$.  It will be convenient 
to refer to a Hilbert space $C$ that has been endowed 
with such a group $W$ as a {\em graded 
Hilbert space}.

In order to discuss free resolutions, we  
shift attention to the more restricted category 
whose objects are graded Hilbert 
modules over the commutative polynomial 
algebra $\mathbb C[z_1,\dots,z_d]$ and whose 
maps are graded morphisms.  In this context, 
one replaces the noncommutative free 
module $F^2=F^2\langle z_1,\dots,z_d\rangle$ with 
its commutative counterpart $H^2=H^2[z_1,\dots,z_d]$, namely 
the completion of $\mathbb C[z_1,\dots,z_d]$ in its 
natural norm.  While this notation differs from 
the notation $H^2(\mathbb C^d)$ used in \cite{arvSubalgIII} 
and \cite{arvCurv}, it is more useful for our purposes here.  
The commutative free Hilbert module $H^2$ is 
realized as a quotient of 
$F^2$ as follows.  
Consider the the operator $U\in\mathcal B(F^2,H^2)$ 
obtained by closing the linear map that sends 
a noncommutative polynomial $f\in \mathbb C\langle z_1,\dots,z_d\rangle$ 
to its commutative image in $\mathbb C[z_d,\dots,z_d]$.  This operator  
is a graded partial isometry whose range is $H^2$ and whose 
kernel is the closure of the commutator ideal in $\mathbb C\langle z_1,\dots,z_d\rangle$, 
\begin{equation}\label{EqB}
K=\overline{{\rm span}}\{f\cdot(z_jz_k-z_kz_j)\cdot g: 
1\leq j,k\leq d,\quad f,g\in \mathbb C\langle z_1,\dots,z_d\rangle\}.  
\end{equation}
One sees this in more concrete terms 
after one identifies $H^2\subseteq F^2$ with the completion of the 
symmetric tensor algebra in the norm inherited from 
$F^2$.  In that realization one has 
$H^2=K^\perp$, and $U$ can be taken as 
the projection with range $K^\perp=H^2$.  
The situation is similar for graded free modules having 
multiplicity; indeed, for any graded coefficient space 
$C$ the map 
$$
U\otimes 1_C: F^2\otimes C\to H^2\otimes C
$$ 
defines a graded cover of the commutative 
free Hilbert module $H^2\otimes C$.

The most general 
graded Hilbert module over the commutative 
polynomial algebra $\mathbb C[z_1,\dots, z_d]$ 
is a graded Hilbert module 
over $\mathbb C\langle z_1,\dots, z_d\rangle$ 
whose underlying row contraction $(T_1,\dots,T_d)$ 
satisfies $T_jT_k=T_kT_j$ for all $j,k$.  
Any vector $\zeta$ in such a module $H$ has a unique decomposition 
into a Fourier series relative 
to the spectral subspaces of the gauge group, 
$$
\zeta=\sum_{n=-\infty}^\infty \zeta_n, 
$$
where $\Gamma(\lambda)\zeta_n=\lambda^n\zeta_n$, 
$n\in\mathbb Z$, $\lambda\in \mathbb T$.   $\zeta$ is said to 
have finite $\Gamma$-spectrum if all but a finite number of 
the terms $\zeta_n$ of this series are zero.  
Finally, a graded contractive module $H$ is 
said to be {\em finitely generated} 
if there is a finite set of vectors 
$\zeta_1,\dots,\zeta_s\in H$, each of 
which has finite $\Gamma$-spectrum, 
such that sums of the form 
$$
f_1\cdot\zeta_1+\cdots+f_s\cdot \zeta_s,
\qquad f_1,\dots,f_s\in\mathbb C\langle z_1,\dots,z_d\rangle
$$
are dense in $H$.

  Our main result is the following 
counterpart of Hilbert's syzygy theorem.  

\begin{thm}\label{thmB}
For every finitely generated 
graded contractive Hilbert module $H$ over the commutative 
polynomial algebra $\mathbb C[z_1,\dots,z_d]$ 
there is a weakly 
exact finite sequence of graded Hilbert modules 
\begin{equation}\label{EqA}
0\longrightarrow F_n\longrightarrow\cdots \longrightarrow 
F_2 \longrightarrow F_1\longrightarrow H\longrightarrow 0    
\end{equation}
in which each $F_k=H^2\otimes C_k$ is a 
free graded commutative Hilbert module of finite 
defect.  
The sequence 
(\ref{EqA}) is unique up to a unitary isomorphism 
of diagrams and it terminates after at most $n=d$ steps.  
\end{thm}

\begin{defn}\label{def2}
The sequence (\ref{EqA}) is called the {\em free} 
resolution of $H$.  
\end{defn}

\begin{rem}[Betti numbers, Euler characteristic]\label{BettiRem}
The sequence (\ref{EqA}) gives rise to a sequence of 
$d$ numerical invariants  of $H$
$$
\beta_k(H)=
\begin{cases}
\defect(F_k),\qquad &1\leq k\leq n,\\
0, \qquad &n<k\leq d.  
\end{cases}
$$
and their alternating sum 
$$
\chi(H)=\sum_{k=1}^d(-1)^{k+1}\beta_k(H)
$$
is called the {\em Euler characteristic} of $H$.  Notice 
that this definition 
makes sense for any finitely generated (graded contractive) 
Hilbert module over $\mathbb C[z_1,\dots,z_d]$, 
and generalizes the 
Euler characteristic of \cite{arvCurv}  
that was 
restricted to Hilbert modules 
of {\em finite rank} as in Remark \ref{rem1}.  
\end{rem}

\begin{rem}[Curvature and Index]
The curvature invariant of \cite{arvCurv} is defined only 
in the context of finite rank contractive Hilbert modules, 
hence the index formula of \cite{arvDirac} that relates the 
curvature invariant to the index of a Dirac operator is 
not meaningful in the broader context of Theorem \ref{thmB}.  
On the other hand, the proof of that formula 
included an argument showing that the Euler characteristic 
can be calculated in terms of the Koszul complex associated 
with the Dirac operator, 
and that part of the proof 
is easily adapted to this context to yield the following 
more general variation of the  
index theorem.  {\em For any finitely generated graded Hilbert 
module  $H$ with Dirac operator $D$, both $\ker D_+$ and 
$\ker D_+^*$ are finite-dimensional, and}
$$
(-1)^d\chi(H)=\dim\ker D_+ - \dim\ker D_+^*.  
$$
\end{rem}

\begin{rem}[Relation to Localized Dilation-Theoretic Resolutions]
We have pointed out in Remark \ref{rem1} 
that for pure $d$-contractions 
$(T_1,\dots,T_d)$ acting on a Hilbert space, dilation-theoretic 
techniques give rise to an exact sequence of contractive 
Hilbert modules and partially isometric maps 
$$
\cdots \longrightarrow H^2\otimes C_2\longrightarrow H^2\otimes C_1
\longrightarrow H\longrightarrow 0, 
$$
in which the coefficient spaces $C_k$ of the free Hilbert modules 
are typically infinite-dimensional, and which apparently 
fails to terminate in a finite number of steps.  However, 
in a recent paper \cite{greene1}, Greene studied 
``localizations" of the above exact sequence at various 
points of the unit ball, and 
he has shown that when one localizes at the origin 
of $\mathbb C^d$, the homology of his localized complex 
agrees with the homology of Taylor's Koszul 
complex (see \cite{taylor1},\cite{taylor2})
of the underlying operator $d$-tuple 
$(T_1,\dots,T_d)$ in all cases.  Interesting as these 
local results are, they appear unrelated to the global 
methods and results of this paper.  
\end{rem}

\begin{rem}[Resolutions of modules over function algebras]
We also point out 
that our use of the terms {\em resolution}  and 
{\em free resolution} differs 
substantially from 
usage of similar terms in work of Douglas, Misra and Varughese 
\cite{dougMV1}, \cite{dougMV2}, \cite{dougMisEQ}, \cite{dougMisQF}.  
For example, in \cite{dougMisQF}, the authors consider 
Hilbert modules over certain algebras $\mathcal A(\Omega)$ 
of analytic functions on 
bounded domains $\Omega\subseteq \mathbb C^d$.  They introduce 
a notion of {\em quasi-free} Hilbert 
module that is related to localization, 
and is characterized as follows.  Consider an inner product 
on the algebraic tensor product $\mathcal A(\Omega)\otimes \ell^2$ 
of vector spaces with three properties: a) evaluations 
at points of $\Omega$ should be locally uniformly bounded, 
b) module multiplication from 
$\mathcal A(\Omega)\times (\mathcal A(\Omega)\otimes\ell^2)$ 
to $\mathcal A(\Omega)\otimes \ell^2$ should be continuous, 
and c) it satisfies a technical condition relating Hilbert norm convergence 
to pointwise convergence throughout $\Omega$.  The completion of 
$\mathcal A(\Omega)\otimes\ell^2$ in that inner product 
gives rise to a Hilbert 
module over $\mathcal A(\Omega)$, and such Hilbert modules 
are called {\em quasi-free}.  

The main result of \cite{dougMisQF} asserts that ``weak" quasi-free 
resolutions 
$$
\cdots\longrightarrow Q_2\longrightarrow Q_1
\longrightarrow H\longrightarrow 0
$$
exist for certain Hilbert modules $H$ over $\mathcal A(\Omega)$, 
namely those that are higher-dimensional generalizations 
of the Hilbert modules studied by 
Cowen and Douglas in \cite{cowDoug} for 
domains $\Omega\subseteq \mathbb C$.  
The modules $Q_k$ are quasi-free 
in the sense above,  
but their ranks may be infinite and 
such sequences may fail to terminate in a finite 
number of steps.  

\end{rem}

\section{Generators}\label{S:Gen}

Throughout this section we consider contractive Hilbert modules 
over the noncommutative 
polynomial algebra $\mathbb C\langle z_1,\dots,z_d\rangle$, 
perhaps graded.  

\begin{defn}
Let $H$ be a Hilbert module over 
$\mathbb C\langle z_1,\dots,z_d\rangle$.  By a {\em generator} 
for $H$ we mean a linear subspace $G\subseteq H$ such that 
$$
H=\overline{\rm span}\,\{f\cdot\zeta: f\in\mathbb C\langle z_1,\dots,z_d\rangle,
\quad \zeta\in G\}.  
$$
\end{defn}

We also say that $H$ is {\em finitely generated} if it has 
a finite-dimensional generator, and 
in the category of graded Hilbert modules 
the term means a bit more, namely, 
that there is a finite-dimensional 
{\em graded} generator.  

According to Definition 
\ref{defStat1}, a finitely generated Hilbert module 
$H$ is properly generated precisely 
when the defect subspace $H\ominus (Z\cdot H)$ is a 
finite-dimensional generator.  
In general, the defect subspace $H\ominus (Z\cdot H)$ of a 
finitely generated Hilbert module is necessarily finite-dimensional, 
but it can fail to generate and is sometimes 
$\{0\}$ (for examples, see the proof 
of Proposition \ref{propGen2}).  
In particular, {\em finitely 
generated Hilbert modules need not be properly generated}.  The 
purpose of this section is to 
show that many important examples 
{\em are} properly generated, 
and that many others are related to properly 
generated Hilbert modules in a simple way.  

The following result can be viewed as 
a noncommutative operator theoretic counterpart of Nakayama's Lemma 
(\cite{EisBk2}, Lemma 1.4).  

\begin{prop}\label{propGen1}
Every finitely generated graded Hilbert module is properly generated.  
\end{prop}

\begin{proof}Let $G=H\ominus(Z\cdot H)$.  $G$ is obviously 
a graded subspace of $H$, and it is finite-dimensional 
because $\dim G=\dim (H/(Z\cdot H))$ is dominated by 
the cardinality of any generating set.  It remains to 
show that $G$ is a generator.

For that, we claim that the 
spectrum of the gauge group $\Gamma$ is bounded below.  
Indeed, the hypothesis implies that 
there is a finite set of elements $\zeta_1,\dots,\zeta_s$ of $H$, 
each having finite $\Gamma$-spectrum, which generate $H$.  By 
enlarging the set of generators appropriately and 
adjusting notation, we can assume 
that each $\zeta_k$ is an eigenvector of $\Gamma$, 
$$
\Gamma(\lambda)\zeta_k=\lambda^{n_k}\zeta_k,\qquad \lambda\in \mathbb T, 
\quad 1\leq k\leq s.  
$$
Let $n_0$ be the minimum of $n_1,n_2,\dots, n_s$.  For any 
monomial $f$ of respective degrees $p_1,\dots, p_d$ in the noncommuting 
variables $z_1,\dots,z_d$ and every $k=1,\dots,s$, $f\cdot \zeta_k$ 
is an eigenvector of $\Gamma$  satisfying 
$$
\Gamma(\lambda)(f\cdot \zeta_k)=\lambda^N f\cdot\zeta_k
$$
with $N=p_1+\cdots+p_d+n_k\geq n_0$.  Since 
elements of this form have $H$ as their closed linear 
span, the spectrum of $\Gamma$ is bounded 
below by $n_0$.  

Setting $H_n=\{\xi\in H: \Gamma(\lambda)\xi=\lambda^n\xi, \lambda\in\mathbb T\}$ 
for $n\in\mathbb Z$, 
we conclude that 
$$
H=H_{n_0}\oplus H_{n_0+1}\oplus\cdots, 
$$
and one has $Z\cdot H_n\subseteq H_{n+1}$ 
for all $n\geq n_0$.  

Since $G=H\ominus(Z\cdot H)$ is gauge-invariant it has a decomposition 
$$
G=G_{n_0}\oplus G_{n_0+1}\oplus \cdots,   
$$
in which $G_{n_0}=H_{n_0}$, $G_n=H_n\ominus (Z\cdot H_{n-1})$ for 
$n>n_0$, and where only a finite number of $G_k$ are 
nonzero.  Thus, each eigenspace $H_n$ decomposes into 
a direct sum 
$$
H_n=G_n\oplus (Z\cdot H_{n-1}),\qquad n>n_0.  
$$
Setting $n=n_0+1$ we have 
$H_{n_0+1}={\rm span}(G_{n_0+1}+Z\cdot G_{n_0})$ and, continuing 
inductively, we find that for all $n>n_0$, 
$$
H_n={\rm span}(G_n+Z\cdot G_{n-1}+Z^{\otimes 2}\cdot G_{n-2}+
\cdots+ Z^{\otimes (n-n_0)}\cdot G_{n_0}),   
$$
where $Z^{\otimes r}$ denotes the space of homogeneous polynomials 
of total degree $r$.  
Since $H$ is spanned by the subspaces $H_n$, 
it follows that $G$ is a generator.  \end{proof}

One obtains the most general examples of 
graded Hilbert submodules of the rank-one free 
commutative Hilbert module $H^2$ in explicit terms 
by choosing a (finite or infinite) 
sequence of homogeneous polynomials $\phi_1,\phi_2, \dots$ 
and considering the closure in $H^2$ of the set of 
all finite linear combinations $f_1\cdot\phi_1 +\cdots+ f_s\cdot\phi_s$, 
where $f_1,\dots,f_s$ are arbitrary polynomials and 
$s=1,2,\dots$.  
In Remark \ref{rem1} above, we alluded 
to the fact that in all nontrivial cases, 
graded submodules of $H^2$ are Hilbert modules of 
infinite rank.  However, 
Proposition \ref{propA} below implies that 
these examples are  properly finitely generated, so 
they have free covers of finite defect by Theorem \ref{thmA}.

\begin{rem}[Examples of Higher Degree]\label{rem2} 
We now describe a class of infinite rank ungraded examples 
with substantially different properties.  Perhaps 
we should point out that there is a more 
general notion of grading with respect to which 
the ambient operators $T_1,\dots,T_d$ in these 
examples are graded 
with various degrees larger than one.  
For brevity, we 
retain the simpler definition 
of grading (\ref{EqStat1}) by viewing 
these examples as ungraded.  
Fix a $d$-tuple of positive 
integers $N_1,\dots,N_d$ and consider 
the following $d$-contraction 
acting on the Hilbert space $H=H^2(\mathbb C^d)$
$$
(T_1, \dots, T_d)=(S_1^{N_1},\dots, S_d^{N_d}),   
$$
where $(S_1,\dots,S_d)$ is the $d$-shift.  The 
defect space of this Hilbert module 
$$
G=H\ominus(T_1H+\cdots+T_dH)^-  
$$
coincides with the intersection of the 
kernels $\ker T_1^*\cap\cdots\cap \ker T_d^*$; 
and in this case one can compute these kernels explicitly,   
with the result 
$$
G={\rm span}\{z_1^{n_1}\cdots z_d^{n_d}: 0\leq n_k < N_k, \quad 1\leq k\leq d\}.  
$$
Moreover, for every set of nonnegative integers $\ell_1,\dots,\ell_d$, 
the set of vectors $T_1^{\ell_1}\cdots T_d^{\ell_d}G$
contains all monomials of the form
$$
z_1^{\ell_1 N_1+n_1}\cdots z_d^{\ell_d N_d+n_d},\qquad 0\leq n_k<N_k, 
\quad 1\leq k\leq d.  
$$
It follows from these observations that $G$ is a proper 
generator for H, and  Theorem \ref{thmA} provides a free cover of the form 
$A: H^2\otimes G\to H$.  

Another straightforward computation with coefficients shows  
that the defect operator of this Hilbert module is of infinite 
rank whenever at least one of the integers $N_1,\dots,N_d$ is 
larger than $1$.  In more detail, each monomial $z^n=z_1^{n_1}\cdots z_d^{n_d}$, 
$n_1,\dots,n_d\geq 0$, is an eigenvector of the defect operator 
$\Delta=(\mathbf 1-T_1T_1^*-\cdots-T_dT_d^*)^{1/2}$; and when $n_k\geq N_k$ for 
all $k$, a straightforward application of the formulas on 
pp. 178--179 of \cite{arvSubalgIII} 
shows that 
$$
\Delta\, z_1^{n_1}\cdots z_d^{n_d}
=c(n)z_1^{n_1}\cdots z_d^{n_d}
$$ 
where the eigenvalues $c(n)=c(n_1,\dots,n_d)$ satisfy $0<c(n)< 1$.   
Hence the defect operator has infinite rank.  
We conclude that, while dilation 
theory provides a coisometry $B: H^2\otimes C\to H$ from 
another free Hilbert module to $H$, it is necessary  
that $C$ be an infinite dimensional Hilbert space.  
Needless to say, 
such a $B$ cannot define a free cover.  
\end{rem}

The preceding examples are all of infinite rank, and it is 
natural to ask about finite rank $d$-contractions -- 
which were the focus of \cite{arvSubalgIII}, 
\cite{arvCurv}, \cite{arvDirac}.  Significantly, while 
the Hilbert module associated with a finite 
rank $d$-contraction is frequently not properly generated, 
it can always be extended to a 
properly generated one by way of a finite-dimensional perturbation.

\begin{prop}\label{propGen2}
Every pure Hilbert module $H$ 
of finite rank can be extended trivially 
to a properly generated one in the sense that there is an 
exact sequence of Hilbert modules 
$$
0\longrightarrow K\longrightarrow \tilde H\underset{A}
\longrightarrow H\longrightarrow 0
$$ 
in which $\tilde H$ is is a properly generated 
pure Hilbert module of the same rank, 
$K$ is a finite-dimensional Hilbert submodule of $\tilde H$, and $A$ is a coisometry.  
\end{prop}

\begin{proof}
A standard dilation-theoretic technique (see \cite{arvSubalgIII} 
for the commutative case, the proof of which works as well in 
general)  
shows that a pure Hilbert module of rank $r$ is 
unitarily equivalent to a quotient of the form 
$$
H\cong (F^2\otimes C)/M
$$
where $F^2$ is the noncommutative free module of 
rank $1$, $C$ is an $r$-dimensional coefficient 
space, and $M$ is a closed submodule of 
$F^2\otimes C$.  We identify $H$ with the orthocomplement 
$M^\perp$ of $M$ in $F^2\otimes C$, with operators 
$T_1,\dots, T_d$ obtained 
by compressing to $M^\perp$ the natural operators 
$U_1\otimes \mathbf 1_C, \dots, U_d\otimes \mathbf 1_C$ 
of $F^2\otimes C$.  

Consider $\tilde H=M^\perp + 1\otimes C$.  This is 
a finite-dimensional extension of $M^\perp$ that is 
also invariant under $U_k^*\otimes \mathbf 1_C$, 
hence it defines a pure Hilbert module of rank $r$ by compressing the 
natural operators in the same way to obtain 
$\tilde T_1,\dots,\tilde T_d\in\mathcal B(\tilde H)$.  Since $\tilde H$ 
contains $H$, the projection $P_{M^\perp}$ restricts to a 
homomorphism of Hilbert modules $A: \tilde H\to H$.   
$A$ is a coisometry, and the kernel of $A$ 
is finite-dimensional because $\dim(\tilde H/H)<\infty$.  

To see that $\tilde H$ is properly generated, one 
computes the defect operator $\Delta$ of $\tilde H$.  Indeed, 
$\Delta=(\mathbf 1_{\tilde H}-\tilde T_1\tilde T_1^*-\cdots-\tilde T_d\tilde T_d^*)^{1/2}$   
is seen to be the compression of the defect operator 
of $U_1\otimes \mathbf 1_C,\dots, U_d\otimes \mathbf 1_C$ 
to $\tilde H$, and the latter defect operator is the 
projection onto $1\otimes C$.  Since $\tilde H$ contains 
$1\otimes C$,  the defect operator of $\tilde H$ is 
the projection on $1\otimes C$.     

Finally, we make use of the observation that a pure finite rank 
$d$-tuple is properly generated whenever its defect operator 
is a projection.  Indeed,  the range of the defect operator 
$\Delta$ 
is always a generator, and when $\Delta$ 
is a projection 
its range coincides with 
$\ker\tilde T_1^*\cap\cdots\cap\tilde \ker T_d^*$.  
\end{proof}

\section{Existence of Free Covers}\label{S:Prf}
We now establish the existence and uniqueness of 
free covers for properly generated 
Hilbert modules over $\mathbb C\langle z_1,\dots,z_d\rangle$.   
A cover $A: F\to H$ 
induces a unitary map of defect spaces; the following 
result implies that this isometry of quotients 
lifts to an isometry of the corresponding subspaces.  

\begin{lem}\label{lemPrf1}
Let $H$ be a Hilbert module and 
let $A: F\to H$ be a cover.  Then 
$A$ restricts to a unitary operator 
from $F\ominus (Z\cdot F)$ to $H\ominus(Z\cdot H)$.  
\end{lem}

\begin{proof}
Let $Q\in\mathcal B(H)$ 
be the projection onto $H\ominus (Z\cdot H)$.  The natural map 
of $H$ onto the quotient Hilbert space $H/(Z\cdot H)$ is a partial isometry 
whose adjoint is the isometry 
$$
 \eta+ Z\cdot H\in H/(Z\cdot H)\mapsto Q\eta\in H\ominus (Z\cdot H), 
\qquad \eta\in H.   
$$
Thus we can define a unitary map  $\tilde A$ from 
$F\ominus(Z\cdot F)$ onto $H\ominus(Z\cdot H)$ by composing 
the three unitary operators 
\begin{align*}
\zeta\in F\ominus (Z\cdot F)&\mapsto \zeta+Z\cdot F\in F/(Z\cdot F),\\
\dot A:F/(Z\cdot F)&\to H/(Z\cdot H),\\
\eta+Z\cdot H\in H/(Z\cdot H)&\mapsto Q\eta\in H\ominus (Z\cdot H), 
\qquad \eta\in H,
\end{align*}
to obtain $\tilde A\zeta=QA\zeta$, $\zeta\in F\ominus (Z\cdot F)$.  
We claim now that $QA\zeta=A\zeta$ for all 
$\zeta\in F\ominus(Z\cdot F)$.  To see that, note that 
$Q^\perp$ is the projection onto $Z\cdot H$, so that for 
all $\zeta\in F\ominus (Z\cdot F)$ one has 
\begin{align*}
\|QA\zeta\|&=\|A\zeta - Q^\perp A\zeta\|=
\inf_{\eta\in Z\cdot H}\|A\zeta-\eta\|\\
&=\|\dot A(\zeta+Z\cdot F)\|_{H/(Z\cdot H)}=
\|\zeta+ Z\cdot F\|_{F/(Z\cdot F)}=\|\zeta\|.  
\end{align*}
Hence, $\|A\zeta-QA\zeta\|^2=\|A\zeta\|^2-\|QA\zeta\|^2=\|A\zeta\|^2-\|\zeta\|^2\leq 0$, 
and the claim follows.  
We conclude that the restriction of $A$ to 
$F\ominus (Z\cdot F)$ is an isometry with range $H\ominus (Z\cdot H)$.  
\end{proof}

\begin{proof}[Proof of Theorem \ref{thmA}]  Let 
$H$ be a properly generated Hilbert module 
over $\mathbb C\langle z_1,\dots,z_d\rangle$ and  
set $C=H\ominus (Z\cdot H)$.  The 
hypothesis asserts that $C$ is a 
generator.  We will show that there is a 
(necessarily unique)  
contraction 
$
A: F^2\otimes C\to H  
$
satisfying 
\begin{equation}\label{EqPrf0}
A(f\otimes \zeta)=f\cdot \zeta,\qquad f\in
\mathbb C\langle z_1,\dots,z_d\rangle , \quad\zeta\in C,  
\end{equation}
and that such an operator $A$ defines 
a free cover.  For that, consider the completely 
positive map defined on $\mathcal B(H)$ by
$\phi(X)=T_1XT_1^*+\cdots+T_dXT_d^*$, and let  
$
\Delta= 
(\mathbf 1-\phi(\mathbf 1))^{1/2}
$
be the defect operator of (\ref{EqStat0}).  Since $H\ominus (Z\cdot H)$ is 
the intersection of kernels $\ker T_1^*\cap\cdots\cap \ker T_d^*=\ker \phi(\mathbf 1)$, 
it follows that 
$$
C=H\ominus (Z\cdot H)=\{\zeta\in H: \Delta\zeta=\zeta\}.  
$$
Thus, $C$ is a 
subspace of the range of $\Delta$ on which 
$\Delta$ restricts to the identity operator, and 
which generates $H$.  
We now use the ``dilation telescope" to show 
that there is a unique contraction 
$L: F^2\otimes \overline{\Delta H}\to H$ such that 
\begin{equation}\label{EqPrf01}
L(f\otimes \zeta)=f\cdot\Delta \zeta,\qquad f\in 
\mathbb C\langle z_1,\dots,z_d\rangle, 
\quad \zeta\in \overline{\Delta H}.    
\end{equation}
Indeed, since 
the monomials 
$\{z_{i_1}\otimes \cdots \otimes z_{i_n}: i_1,\dots,i_n\in\{1,\dots,d\}\}$, 
$n=1,2,\dots$, 
together with the constant polynomial $1$, 
form an orthonormal basis for $F^2$, 
the formal adjoint of $L$ is easily computed and found to be 
$$
L^*\xi = 1\otimes \Delta\xi + \sum_{n=1}^\infty\sum_{i_1,\dots,i_n=1}^d
z_{i_1}\otimes \cdots\otimes z_{i_n}\otimes \Delta T_{i_n}^*\cdots T_{i_1}^*\xi, 
\qquad \xi\in H.  
$$
One calculates norms in the obvious way to obtain
\begin{align*}
\|L^*\xi\|^2 &=  \|\Delta\xi\|^2+\sum_{n=1}^\infty\sum_{i_1,\dots,i_n=1}^d
\|\Delta T_{i_n}^*\cdots T_{i_1}^*\xi\|^2\\
&=\langle (\mathbf 1-\phi(\mathbf 1))\xi,\xi\rangle +
\sum_{n=1}^\infty \langle (\phi^n(\mathbf 1-\phi(\mathbf 1))\xi,\xi\rangle
\\&=\langle (\mathbf 1-\phi(\mathbf 1))\xi,\xi\rangle +
\sum_{n=1}^\infty \langle (\phi^n(\mathbf 1)-\phi^{n+1}(\mathbf 1))\xi,\xi\rangle
\\&=\|\xi\|^2-\lim_{n\to\infty}\langle \phi^n(\mathbf 1)\xi,\xi\rangle\leq \|\xi\|^2.   
\end{align*}
Hence $\|L\|=\|L^*\|\leq 1$.   
Finally, let $A$ be the restriction of $L$ to 
the submodule 
$F^2\otimes C\subseteq F^2\otimes \overline{\Delta H}$, 
where we now consider $F^2\otimes C$ as a free Hilbert module 
of possibly smaller defect.  
Since $\Delta$ restricts to the identity 
on $C$, (\ref{EqPrf0}) follows from (\ref{EqPrf01}).  

By its definition, the restriction of $A$ to $1\otimes C$ is 
an isometry with range $C=H\ominus (Z\cdot H)$, 
hence $A$ induces a unitary operator 
of defect spaces 
$$
\dot A: (F^2\otimes C)/(Z\cdot(F^2\otimes C))\cong1\otimes C 
\to C= H\ominus (Z\cdot H)\cong H/(Z\cdot H).  
$$
The range of $A$ is dense, since it contains 
$$
\{f\cdot \zeta: f\in\mathbb C\langle z_1,\dots,z_d\rangle, \ \zeta\in H\ominus (Z\cdot H)\}
$$
and $H$ is properly generated.  
It follows that $A:F^2\otimes C\to H$ is a free cover.  

For uniqueness, 
let $B: \tilde F= F^2\otimes \tilde C\to H$ be another free cover 
of $H$.  We exhibit a unitary isomorphism of Hilbert modules 
$V\in\mathcal B(F^2\otimes \tilde C,F^2\otimes C)$ such 
that $BV=A$ as follows.  We have already pointed out that the defect space 
of $\tilde F=F^2\otimes \tilde C$ (resp. $F=F^2\otimes C)$) 
is identified with $1\otimes \tilde C$ 
(resp. $1\otimes C$).  Similarly, the defect space of $H$ is 
identified with $H\ominus (Z\cdot H)$.  
Since both $A$ and $B$ are  
covers of $H$, Lemma \ref{lemPrf1} implies 
that they restrict to unitary operators, from 
the respective spaces $1\otimes C$ and $1\otimes\tilde C$, onto the same  
subspace $H\ominus (Z\cdot H)$ of $H$.  Thus 
there is a unique unitary operator $V_0: C\to \tilde C$ 
that satisfies 
$$
A(1\otimes\zeta)=B(1\otimes V_0\zeta),\qquad \zeta\in C.  
$$
Let $V=\mathbf 1_{F^2}\otimes V_0\in
\mathcal B(F^2\otimes\tilde C,F^2\otimes C)$.  Obviously 
$V$ is a unitary operator, and it satisfies $BV=A$ since 
for every polynomial $f\in\mathbb C\langle z_1,\dots,z_d\rangle$ 
$$
BV(f\otimes\zeta)=B(f\cdot (1\otimes V_0\zeta))=
f\cdot B(1\otimes V_0\zeta)=f\cdot A(1\otimes \zeta)=A(f\otimes\zeta), 
$$
and one can take the closed linear span 
on both sides.  
$V$ must implement an isomorphism 
of modules since for any polynomials 
$f,g\in\mathbb C\langle z_1,\dots,z_d\rangle$ 
and every $\zeta\in \tilde C$ we have 
$$
V(f\cdot(g\otimes \zeta))=(\mathbf 1\otimes V_0)(f\cdot g\otimes\zeta)
=f\cdot (g\otimes V_0\zeta)=f\cdot V(g\otimes \zeta).   
$$

Conversely, if a Hilbert module $H$ has a free 
cover $A: F^2\otimes C\to H$, then since 
$1\otimes C$ is the orthocomplement of 
$Z\cdot(F^2\otimes C)$, Lemma \ref{lemPrf1} implies 
that $A(1\otimes C)=H\ominus (Z\cdot H))$.  Since 
$A$ is a module homomorphism, we see that 
$$
A({\rm span}\{f\otimes \zeta: f\in\mathbb C,\quad \zeta\in C\})
={\rm span}\{f\cdot \zeta: f\in\mathbb C,\quad \zeta\in H\ominus (Z\cdot H)\}.
$$  
The closure of the left side is $H$ because $A$ has dense 
range, and we conclude that $H\ominus (Z\cdot H)$ is a generator
of $H$.  
\end{proof}

We require the following 
consequence of Theorem \ref{thmA} for finitely 
generated graded Hilbert modules.  

\begin{thm}\label{thmC}
Every finitely generated graded Hilbert module $H$ over 
the noncommutative polynomial algebra 
has a graded free cover $A: F^2\otimes C\to H$, 
and any two graded graded free covers are equivalent.  

If the underlying operators of $H$ commute, then this  
free cover descends naturally to a commutative 
graded free cover $B: H^2\otimes C\to H$.  
\end{thm}

\begin{proof} 
Proposition \ref{propGen1} implies that the 
space $C=H\ominus (Z\cdot H)$ is a 
finite-dimensional generator.  Moreover, 
since $Z\cdot H$ is invariant under the gauge group, 
so is $C$, and
the restriction of the gauge 
group to $C$ gives rise to a unitary 
representation $W: \mathbb T\to \mathcal B(C)$ of the 
circle group on $C$.  

If we make the free Hilbert module $F^2\otimes C$ into a graded 
one by introducing the gauge group 
$$
\Gamma(\lambda)=\Gamma_0(\lambda)\otimes W(\lambda),
\qquad \lambda \in \mathbb T, 
$$
then we claim that 
the map $A: F^2\otimes C\to H$ defined in the proof 
of Theorem \ref{thmA} must intertwine 
$\Gamma$ and $\Gamma_H$.  Indeed,  this follows 
from the fact that for 
every polynomial $f\in\mathbb C\langle z_1,\dots,z_d\rangle$,
every   
$\zeta\in C$, and every $\lambda\in\mathbb T$, one has 
\begin{align*}
\Gamma_H(\lambda)A(f\otimes \zeta)&=
\Gamma_H(\lambda)(f\cdot\zeta)=
f(\lambda z_1,\dots,\lambda z_d)\cdot \Gamma_H(\lambda)\zeta\\&=
A(\Gamma_0(\lambda)f\otimes W(\lambda)\zeta)=
A\Gamma(\lambda)(f\otimes\zeta).
\end{align*}

The proof of uniqueness in the graded context 
is now a straightforward 
variation of the uniqueness proof of 
Theorem \ref{thmA}.  Finally, since 
$H^2$ is naturally identified with the quotient 
$F^2/K$ where $K$ is the closure of the commutator 
ideal in $\mathbb C\langle z_1,\dots,z_d\rangle$, 
it follows that when the underlying operators commute, 
the cover $A: F^2\otimes C\to H$ factors naturally through 
$(F^2/K)\otimes C \sim H^2\otimes C$ and one 
can promote $A$ to a graded 
commutative free cover $B: H^2\otimes C\to H$.  
\end{proof}

\section{Existence of Free Resolutions}\label{S:res}

We turn now to the proof of existence of finite resolutions for 
graded Hilbert modules over the 
commutative polynomial algebra $\mathbb C[z_1,\dots,z_d]$.   
We require some algebraic results obtained by 
Hilbert at the end of the century before 
last \cite{hilb1}, \cite{hilb2}.  While Hilbert's theorems have 
been extensively generalized, what we require are the 
most concrete 
versions of a) the basis theorem and b) the syzygy theorem.  
We now describe these 
classical results in a formulation that is convenient for our 
purposes, referring 
the reader to \cite{northFree},  \cite{EisBk} and \cite{serLoc} for 
more detail on the underlying linear algebra.  

Let $T_1,\dots, T_d$ be a set of commuting linear operators 
acting on a complex vector space $M$.  We view $M$ as a module 
over $\mathbb C[z_1,\dots,z_d]$ in the usual way, with 
$f\cdot \xi=f(T_1,\dots,T_d)\xi$, $f\in\mathbb C[z_1,\dots,z_d]$, 
$\xi\in M$.  Such a module is said to be {\em graded} if there 
is a specified sequence $M_n$, $n\in\mathbb Z$, of 
subspaces that gives rise to an algebraic direct sum decomposition 
$$
M=\sum_{n=-\infty}^\infty M_n
$$
with the property $T_kM_n\subseteq M_{n+1}$, for all $k=1,\dots,d$, 
$n\in\mathbb Z$.  Thus, every element $\xi$ of $M$ admits a unique 
decomposition $\xi=\sum_n\xi_n$, where $\xi_n$ belongs to $M_n$ and 
$\xi_n=0$ for all but a finite number of $n$.  
We confine ourselves to the standard 
grading on $\mathbb C[z_1,\dots,z_d]$ in which the generators 
$z_1,\dots,z_d$ are all of degree $1$.  Finally, $M$ is said to be 
{\em finitely generated} if there is a finite set $\zeta_1,\dots,\zeta_s\in M$ 
such that 
$$
M=\{f_1\cdot \zeta_1+\cdots+f_s\cdot \zeta_s: 
f_1,\dots,f_d\in\mathbb C[z_1,\dots,z_d]\}.  
$$
A free module is a module of the 
form $F=\mathbb C[z_1,\dots,z_d]\otimes C$ 
where $C$ is a complex vector space, 
the module action being defined in the usual way by 
$f\cdot(g\otimes \zeta)=(f\cdot g)\otimes \zeta$.  The 
{\em rank} of $F$ is the dimension of $C$.  
A free module 
can be graded in many ways, and for our purposes 
the most general grading on 
$F=\mathbb C[z_1,\dots,z_d]\otimes C$ is defined 
as follows.  Given an arbitrary grading on the ``coefficient" 
vector space $C$ 
$$
C=\sum_{n=-\infty}^\infty C_n, 
$$
there is a corresponding grading of the tensor product 
$F=\sum_n F_n$ in which 
$$
F_n=\sum_{k=0}^\infty Z^k\otimes C_{n-k},   
$$
where $Z^k$ denotes the space of all homogeneous 
polynomials of degree $k$ in $\mathbb C[z_1,\dots,z_d]$, 
and where the sum on the right denotes the linear subspace of $F$ 
spanned by $\cup\{Z^k\otimes C_{n-k}: k\in \mathbb Z\}$.  
If $C$ is finite-dimensional, then there are integers 
$n_1\leq n_2$ such that 
$$
C=C_{n_1}+C_{n_1+1}+\cdots+ C_{n_2}, 
$$
so that 
\begin{equation}\label{EqPrf1}
F_n=\sum_{k=0}^\infty Z^k\otimes C_{n-k}=
\sum_{k=\max(n-n_2,0)}^{\max(n-n_1,0)} Z^k\otimes C_{n-k}
\end{equation}
is finite-dimensional for each $n\in\mathbb Z$, 
$F_n=\{0\}$ for $n<n_1$, and 
$F_n$ is spanned by $Z^{n-n_2}\cdot F_{n_2}$ for $n\geq n_2$.  

Homomorphisms 
of graded modules $u: M\to N$ are required to be of degree zero 
$$
u(M_n)\subseteq N_n,\qquad n\in\mathbb Z.  
$$
It will also be convenient to adapt 
Serre's definition of minimality for 
homomorphisms of modules over 
local rings (page 84 of \cite{serLoc}) to homomorphisms of graded 
modules over $\mathbb C[z_1,\dots,z_d]$, as follows.   
A homomorphism $u: M\to N$ of modules is said 
to be {\em minimal} if it induces an isomorphism 
of vector spaces 
$$
\dot u: M/(z_1\cdot M+\cdots+z_d\cdot M)\to 
u(M)/u(z_1\cdot M+\cdots+z_d\cdot M).
$$
Equivalently, $u$ is minimal iff 
$\ker u\subseteq z_1\cdot M+\cdots+z_d\cdot M$.

A {\em free resolution} of an algebraic 
graded module $M$ is a (perhaps infinite) 
exact sequence of graded modules 
$$
\cdots\longrightarrow F_n\longrightarrow \cdots\longrightarrow 
F_2\longrightarrow F_1\longrightarrow M\longrightarrow 0, 
$$
where each $F_r$ is free or $0$.  Such a resolution 
is said to be {\em finite} 
if each $F_r$ is of finite rank and  
$F_r=0$ for sufficiently large $r$, and {\em minimal} if 
for every $r=1,2,\dots$, the arrow emanating from  
$F_r$ denotes a minimal homomorphism.  

\begin{thm}[Basis Theorem]\label{basisThm}
Every submodule of a finitely generated module over 
$\mathbb C[z_1,\dots,z_d]$ 
is finitely generated.  
\end{thm}

\begin{thm}[Syzygy Theorem]\label{syzThm}
Every finitely generated graded module $M$ over $\mathbb C[z_1,\dots,z_d]$ 
has a finite free resolution 
$$
0\longrightarrow F_n\longrightarrow \cdots\longrightarrow 
F_2\longrightarrow F_1\longrightarrow M\longrightarrow 0  
$$
that is minimal with length $n$ at most $d$, and any two 
minimal resolutions are isomorphic.  
\end{thm}

While we have stated the ungraded version of the basis theorem, 
all we require is the special case for graded modules.  
We
base the proof of Theorem \ref{thmB} on two operator-theoretic 
results, the first of which is a Hilbert space counterpart of 
the basis theorem for graded modules.  

\begin{prop}\label{propA}
Let $H$ be a finitely generated graded Hilbert module 
over $\mathbb C[z_1,\dots,z_d]$ and let 
$K\subseteq H$ be a 
closed gauge-invariant submodule.  Then $K$ 
is a properly 
generated graded Hilbert module 
of finite defect.  
\end{prop}

\begin{proof}We first collect some structural information 
about $H$ itself.  Let $\Gamma$ be the gauge group 
of $H$ and consider the spectral subspaces of $\Gamma$ 
$$
H_n=\{\xi\in H: \Gamma(\lambda)\xi=\lambda^n\xi\},\qquad n\in \mathbb Z.  
$$
The finite-dimensional subspace 
$G=H\ominus (Z\cdot H)$ is invariant 
under the action of $\Gamma$, and Proposition 
\ref{propGen1} implies that $G$ is a generator.   Writing 
$G_n=G\cap H_n$, $n\in\mathbb Z$, it follows that $G$ decomposes 
into a finite sum of mutually orthogonal subspaces 
$$
G=G_{n_1}+G_{n_1+1}+\cdots+G_{n_2},  
$$
where $n_1\leq n_2$ are fixed integers.  
A computation similar to that of (\ref{EqPrf1}) shows 
that $H_n=\{0\}$ for $n<n_1$, and for 
$n\geq n_1$, $H_n$ can be expressed in terms of the $G_k$ by way of 
\begin{equation}\label{EqPrf2}
H_n=\sum_{k=\max(n-n_2,0)}^{\max(n-n_1,0)} Z^k\cdot G_{n-k}. \qquad n\in\mathbb Z;
\end{equation}
in particular, each $H_n$ is finite-dimensional.  

Consider the algebraic module 
$$
H^0={\rm span}\{f\cdot\zeta: f\in\mathbb C[z_1,\dots,z_d],\quad\zeta\in G\}
$$
generated by $G$.  Formula (\ref{EqPrf2}) shows that $H^0$ is 
linearly spanned by the spectral subspaces of $H$, 
\begin{equation}\label{EqPrf3}
H^0=H_{n_1}+H_{n_1+1}+\cdots.  
\end{equation}

Now let $K\subseteq H$ be a closed invariant subspace that is also 
invariant under the action of $\Gamma$.  Letting $K_n=H_n\cap K$ be 
the corresponding spectral subspace of $K$, then we have a decomposition of 
$K$ into mutually orthogonal finite-dimensional subspaces 
$$
K=K_{n_1}\oplus K_{n_1+1}\oplus \cdots,   
$$
such that $z_kK_n\subseteq K_{n+1}$, for $1\leq k\leq d$, $n\geq n_1$.  
Let $K^0$ be the (nonclosed) linear span 
$$
K^0=K_{n_1}+ K_{n_1+1}+ \cdots.
$$
Obviously, $K^0$ is dense in $K$ and it is a submodule of 
the finitely generated algebraic module $H^0$.  
Theorem 
\ref{basisThm} implies that there is a finite set of vectors 
$\zeta_1,\dots,\zeta_s\in K^0$ such that 
$$
K^0=\{f_1\cdot \zeta_1+\cdots+f_s\cdot \zeta_s: 
f_1,\dots,f_s\in\mathbb C[z_1,\dots,z_d]\}.  
$$
Choosing $p$ large enough that $\zeta_1,\dots,\zeta_s\in K_{n_1}+\cdots+K_p$, 
we find that $K_{n_1}+\cdots +K_p$ is a graded finite-dimensional generator 
for $K$.  An application of Proposition \ref{propGen1} now completes 
the proof.  
\end{proof}

\subsection{From Hilbert Modules to Algebraic Modules}

A finitely generated 
graded Hilbert module $H$ over $\mathbb C[z_1,\dots,z_d]$ has 
many finite-dimensional graded generators $G$; 
if one fixes such a $G$ then 
there is an associated algebraic graded module 
$M(H,G)$ over $\mathbb C[z_1,\dots,z_d]$, namely 
$$
M(H,G)={\rm span}\{f\cdot \zeta: f\in\mathbb C[z_1,\dots,z_d], \quad \zeta\in G\}.
$$  
The second 
result that we require is that it is possible to make 
appropriate choices of $G$ so as to obtain a functor from Hilbert 
modules to algebraic modules.  We now define this functor and collect 
its basic properties.  

Consider the category $\mathcal H_d$ 
whose objects are graded finitely generated 
Hilbert modules over $\mathbb C[z_1,\dots,z_d]$, 
with covers as maps.  Thus, 
$\hom(H,K)$ consists of 
graded homomorphisms $A: H\to K$ satisfying $\|A\|\leq 1$,  
such that $AH$ is dense in $K$, and which induce unitary operators of 
defect spaces 
$$
\dot A: H/(Z\cdot H)\to K/(Z\cdot K).
$$
Since we are requiring maps in $\hom(H,K)$ to have dense 
range, a straightforward argument (that we omit) shows 
that $\hom(\cdot,\cdot)$ is 
closed under composition.  

The corresponding algebraic category $\mathcal A_d$ 
has objects consisting of 
graded finitely generated modules over $\mathbb C[z_1,\dots,z_d]$,  
in which $u\in \hom(M,N)$ means that $u$ is a {\em minimal} 
graded homomorphism satisfying $u(M)=N$.

\begin{prop}\label{lemB}
For every Hilbert module $H$ in $\mathcal H_d$ let $H^0$ be the 
algebraic module over $\mathbb C[z_1,\dots,z_d]$ defined by 
$$
H^0={\rm span}\{f\cdot \zeta: f\in\mathbb C[z_1,\dots,z_d],\quad 
\zeta\in H\ominus (Z\cdot H)\}.
$$
Then $H^0$ belongs to $\mathcal A_d$.  Moreover, 
for every $A\in \hom(H,K)$ one has $AH^0= K^0$, 
and the restriction $A^0$ of $A$ to $H^0$ defines an element 
of $\hom(H^0,K^0)$.  
The association $H\to H^0$, $A\to A^0$ is a functor 
satisfying: 
\begin{enumerate}
\item[(i)] For every $H\in \mathcal H_d$, $H^0=\{0\}\implies H=\{0\}$.  
\item[(ii)]For every $A\in \hom(H,K)$, $A^0=0\implies A=0$.  
\item[(iii)] For every free graded Hilbert module $F$ of defect $r$, 
$F^0$ is a free algebraic graded module of rank $r$.  
\end{enumerate}
\end{prop}

\begin{proof}
Since the defect subspace $H\ominus (Z\cdot H)$ is finite-dimensional 
and invariant under the action of the gauge group $\Gamma$, $H^0$ is a 
finitely generated module over the polynomial algebra that is 
invariant under the action of the gauge group.  Thus it acquires 
an algebraic grading $H^0=\sum_n H^0_n$ by setting 
$$
H^0_n=H^0\cap H_n=\{\xi\in H^0: \Gamma(\lambda)\xi=\lambda^n\xi, 
\quad \lambda\in\mathbb T\}, \qquad n\in\mathbb Z.  
$$

Let $H,K\in\mathcal H_d$ and let $A\in\hom(H,K)$.  
Lemma \ref{lemPrf1} implies that 
$$
A(H\ominus (Z\cdot H))=K\ominus (Z\cdot K),
$$ 
so 
that $A$ restricts to a surjective graded homomorphism of 
modules $A^0: H^0\to K^0$.  We claim that 
$A^0$ is minimal, i.e., $\ker A^0\subseteq z_1\cdot H^0+\cdots+z_d\cdot H^0$.  
To see that, choose $\xi\in H^0$ such that $A\xi=0$.  Since $H^0$ decomposes into a sum
$$
H^0=H\ominus (Z\cdot H) + z_1\cdot H^0+\cdots +z_d\cdot H^0,   
$$
we can decompose $\xi$ correspondingly 
$$
\xi = \zeta + z_1\cdot \eta_1+\cdots +z_d\cdot \eta_d, 
$$
where $\zeta\in H\ominus (Z\cdot H)$ and $\eta_j\in H^0$.  
Since $\dot A$ is an injective operator defined on $H/(Z\cdot H)$, 
$\ker A$ must be contained in $Z\cdot H$.  It follows that $\xi\in Z\cdot H$, 
and therefore $\zeta=\xi-z_1\cdot\eta_1-\cdots-z_d\cdot\eta_d\in 
Z\cdot H=(H\ominus (Z\cdot H))^\perp$ 
is orthogonal to itself.  Hence $\zeta=0$, and 
we have the desired conclusion 
$$
\xi=z_1\cdot \eta_1+\cdots+z_d\cdot\eta_k\in z_1\cdot H^0+\cdots +z_d\cdot H^0. 
$$ 
The restriction $A^0$ of $A$ to $H^0$ is therefore a minimal 
homomorphism, whence $A^0\in\hom(H^0,K^0)$.   

The composition rule $(AB)^0=A^0B^0$ follows immediately, so that 
we have defined a functor.  Finally, 
both properties (i) and (ii) are consequences of the fact that, by 
Proposition \ref{propGen1}, $H^0$ is dense in $H$, while  
(iii) is obvious.  
\end{proof}

\begin{proof}[Proof of Theorem \ref{thmB}]
Given a graded 
finitely generated Hilbert module $H$, we claim 
that there is a weakly exact sequence 
\begin{equation}\label{EqPrf4}
\cdots \longrightarrow F_n\longrightarrow \cdots \longrightarrow 
F_2\longrightarrow F_1\longrightarrow H\longrightarrow 0,
\end{equation}
in which each $F_r$ is a free graded Hilbert module of finite defect.  
Indeed, Proposition \ref{propA} implies that $H$ is properly 
generated, and  by Theorem \ref{thmB}, it has a graded free cover
$A: F_1\to H$ in which $F_1=H^2\otimes C_1$ is a graded 
free Hilbert module with 
$\dim C_1=\defect(F_1)=\defect(H)<\infty$.  
This gives a sequence of graded Hilbert modules 
\begin{equation}\label{EqPrf5}
F_1\underset{A}\longrightarrow H\longrightarrow  0
\end{equation}
that is weakly exact at $H$.  Proposition \ref{propA} implies 
that $\ker A$ is a 
properly generated graded Hilbert module of finite defect, 
so that another application of Theorem \ref{thmB} produces 
a graded free cover $B: F_2\to \ker A$ in which $F_2$ is 
a graded free Hilbert module of finite defect.  
Thus we can extend 
(\ref{EqPrf5}) to a longer sequence 
$$
F_2\longrightarrow F_1\longrightarrow H\longrightarrow 0
$$
that is weakly exact at $F_1$ and $H$.  Continuing 
inductively, we obtain (\ref{EqPrf4}).  

Another application of 
Theorem \ref{thmB} implies that the sequence 
(\ref{EqPrf4}) is uniquely determined by $H$ up to a 
unitary isomorphism of diagrams.  
The only issue remaining is whether its length is finite.  
To see that (\ref{EqPrf4}) must terminate, 
consider the associated sequence 
of graded algebraic modules provided by Proposition \ref{lemB}
$$
\cdots \longrightarrow F_n^0\longrightarrow \cdots \longrightarrow 
F_2^0\longrightarrow F_1^0\longrightarrow H^0\longrightarrow 0.  
$$
Proposition \ref{lemB} implies that 
this is a {\em minimal} free resolution of $H^0$  into 
graded free modules $F_r^0$ of finite rank.  The 
uniqueness assertion of Theorem \ref{syzThm} 
implies that there is an integer 
$n\leq d$ such that $F_r^0=0$ for all $r>n$.  By 
Proposition \ref{lemB} (i), we have $F_r=0$ for $r>n$.  
\end{proof}

\begin{rem}[Noncommutative Generalizations]
Perhaps it is worth pointing out that there is 
no possibility of generalizing Theorem \ref{thmB} to 
the noncommutative setting, the root cause being that Hilbert's 
basis theorem fails for modules over 
the noncommutative algebra 
$\mathbb C\langle z_1,\dots,z_d\rangle$.  More precisely, 
there are finitely generated graded Hilbert 
modules $H$ over 
$\mathbb C\langle z_1,\dots,z_d\rangle$ that 
do not have finite free resolutions.  Indeed, 
while Theorem \ref{thmA} implies that for any 
such Hilbert module $H$ there is 
a graded free Hilbert module 
$F_1=F^2\otimes C$ with $\dim C<\infty$ 
and a weakly exact sequence 
of graded Hilbert modules 
$$
F_1\underset{A}\longrightarrow H\longrightarrow 0, 
$$
and while the kernel of $A$ is a certainly a 
graded submodule of 
$F^2\otimes C$, the kernel of $A$ need not 
be finitely generated.  For such a Hilbert module $H$, this 
sequence cannot 
be continued beyond $F_1$ within the category of Hilbert modules 
of finite defect.  

As a concrete example of this phenomenon, let $N\geq 2$ be an integer,  
let $Z=\mathbb C^d$ for some $d\geq 2$, and consider 
the free graded noncommutative Hilbert module 
$$
F^2=\mathbb C\oplus Z\oplus Z^{\otimes 2}\oplus 
Z^{\otimes 3}\oplus \cdots.
$$  
We claim that there is an infinite 
sequence of unit vectors $\zeta_N, \zeta_{N+1},\cdots\in F^2$ such 
that $\zeta_n\in Z^{\otimes n}$ and, for all $n\geq N$, 
$$
\zeta_{n+1}\perp M_n=\{f_N\cdot\zeta_N+f_{N+1}\cdot\zeta_{N+1}+\cdots+f_n\cdot\zeta_{n}: 
f_N,\dots,f_n\in\mathbb C\langle z_1,\dots,z_d\rangle\}.
$$  
Indeed, choose a unit vector $\zeta_N$ arbitrarily in $Z^{\otimes N}$ and, 
assuming that $\zeta_N,\dots, \zeta_n$ have been defined with the 
stated properties, 
note that $M_n$ is a graded submodule of $F^2$ 
such that 
$$
M_n\cap Z^{\otimes(n+1)}=
Z^{\otimes(n+1-N)}\cdot \zeta_N+Z^{\otimes(n-N)}\cdot \zeta_{N+1}
+\cdots+Z\cdot \zeta_{n}. 
$$
Recalling that $\dim Z^k=d^k$, an obvious dimension estimate implies that 
\begin{align*}
\dim(M_n\cap Z^{\otimes(n+1)})&\leq d^{n+1-N}+\cdots+d=d\frac{d^{n-N+1}-1}{d-1}\\
&< \frac{d^{n-N+2}}{d-1}\leq  d^{n-N+2}< d^{n+1}
=\dim(Z^{\otimes (n+1)}).  
\end{align*}
Hence there is a unit vector $\zeta_{n+1}\in Z^{\otimes(n+1)}$ 
that is orthogonal to $M_n$.  
Now let $M$ be the closure of $M_N\cup M_{N+1}\cup\cdots$.  $M$ is a 
graded invariant subspace of $F^2$ with the property that 
$M\ominus (Z\cdot M)$ contains the orthonormal set 
$\zeta_N, \zeta_{N+1},\dots$,  so that $M$ cannot be 
finitely generated.  

Finally, if we take $H$ to be the 
Hilbert space quotient $F^2/M$, then $H$ is a graded 
Hilbert module over 
$\mathbb C\langle z_1,\dots,z_d\rangle$ having a single 
gauge-invariant 
cyclic vector $1+M$, such that the 
natural projection $A: F^2\to H=F^2/M$ is a graded 
free cover of $H$ where 
$\ker A=M$ is not finitely generated.  
\end{rem}

Notice that the preceding construction used the fact that the 
dimensions of the spaces 
$Z^{\otimes k}$ of 
noncommutative homogeneous polynomials 
grow exponentially in $k$.  If one attempts 
to carry out this construction in the commutative setting,  
in which $F^2$ is replaced by $H^2$, one will find that 
the construction of the sequence $\zeta_N,\zeta_{N+1},\dots$ 
fails at some point because the dimensions of the 
spaces $Z^k$ of homogeneous polynomials 
grow too slowly.  Indeed, as reformulated 
in Proposition \ref{propA}, 
Hilbert's remarkable basis theorem implies that this 
construction  {\em must} fail in the commutative 
setting, since every graded 
submodule of $H^2$ is finitely generated.  

\section{Examples of Free Resolutions}\label{S:examp}

In this section we discuss some examples 
of free resolutions and their 
associated Betti numbers.  
There are two simple - and 
closely related - procedures for 
converting a free Hilbert module into one that is 
no longer free, by changing its ambient operators 
as follows.  
\begin{enumerate}
\item Append a number $r$ of zero operators to 
the $d$-shift 
$(S_1,\dots,S_d)$ to obtain 
a $(d+r)$-contraction acting on 
$H^2[z_1,\dots,z_d]$ that is not the 
$(d+r)$-shift.
\item Pass 
from $H^2[z_1,\dots,z_d]$ to a 
quotient $H^2[z_1,\dots,z_d]/K$ 
where $K$ is the closed submodule generated by 
some of the coordinates $z_1,\dots, z_d$.  
\end{enumerate}

We begin by pointing out that one can 
understand either of these examples 
(1) or (2) by analyzing 
the other.  We then calculate the Betti numbers of 
the Hilbert modules of (1) in the case where 
one appends three zero operators to the $d$-shift.  
In order to calculate the Betti numbers 
of a graded Hilbert module one has to calculate 
its free resolution, and that is the route we 
follow.  

To see that (1) and (2) are equivalent 
constructions, consider the operator 
$(d+r)$-tuple $\bar T=(S_1,\dots, S_d, 0, \dots,0)$
obtained from the $d$-shift 
$(S_1,\dots,S_{d})$ acting 
on $H^2[z_1,\dots,z_{d}]$ 
by adjoining $r$ zero operators.   
Let $K$ 
be the closed invariant subspace 
of $H^2[z_1,\dots,z_{d+r}]$ generated by 
$z_{d+1},z_{d+2},\dots,z_{d+r}$.
Recalling that 
$H^2[z_1,\dots,z_d]$ embeds isometrically  
in $H^2[z_1,\dots, z_{d+r}]$ with orthocomplement $K$,
$$
H^2[z_1,\dots, z_{d+r}]=H^2[z_1,\dots, z_d]\oplus K,   
$$
one finds that the quotient Hilbert module 
$H^2[z_1,\dots,z_{d+r}]/K$ is identified with 
$H^2[z_1,\dots,z_d]$ in such a way that the 
natural $(d+r)$-contraction defined by 
the quotient is unitarily 
equivalent to $\bar T$. 

Before turning to explicit computations we emphasize 
that, in order to calculate free resolutions, 
one has to iterate the process of calculating 
free covers.  We begin by summarizing 
that procedure in explicit terms.  

\begin{rem}[Free Covers and Free Resolutions]\label{rem2}
Let $H$ be a finitely generated graded Hilbert module 
over $\mathbb C[z_1,\dots,z_d]$.  In order to calculate the 
free resolution of $H$ one has to iterate the following 
procedure.  

\begin{enumerate}
\item
One first calculates the free cover
$A_1: H^2[z_1,\dots,z_d]\otimes G_1\to H$  of 
the given Hilbert module $H$, following the 
proof of Theorem \ref{thmA}.  
More explicitly, one calculates 
the unique proper generator 
$G_1\subseteq H$ 
$$
G_1=H\ominus(Z\cdot H),   
$$
the connecting map $A_1$ being the closure of  
the multiplication map 
$$
A(f\otimes \zeta)=f\cdot\zeta,\qquad f\in \mathbb C[z_1,\dots, z_d],
\quad \zeta\in G_1,   
$$
where the free Hilbert module 
$H^2[z_1,\dots, z_d]\otimes G_1$ is endowed with the 
grading $\Gamma_0\otimes W$,  $W$ being the 
unitary representation of the circle group on $G$ defined 
by restricting the grading $\Gamma_H$ of $H$,
$$
W(\lambda)=\Gamma_H(\lambda)\restriction_G,\qquad \lambda\in\mathbb T.  
$$
Notice that in order to carry out this step, one 
simply has to identify 
$Z\cdot H$ and its orthocomplement in concrete terms.  
\item
One then replaces $H$ with the finitely generated graded Hilbert 
module $\ker A_1\subseteq H^2[z_1,\dots,z_d]\otimes G_1$ and 
repeats the procedure.  
It is significant that in order to 
continue, one must 
identify the kernel of $A_1$ 
and its proper generator 
$G_2=\ker A_1\ominus (Z\cdot \ker A_1)$.  
\end{enumerate}
According to Theorems 
\ref{thmA} and \ref{thmB}, this process will  
terminate in the zero Hilbert module 
after at most $d$ steps, and the resulting 
sequence 
$$
0\longrightarrow H^2[z_1,\dots,z_d]\otimes G_n
\underset{A_n}\longrightarrow
\cdots\underset{A_2}\longrightarrow H^2[z_1,\dots,z_d]\otimes G_1
\underset{A_1}\longrightarrow H\longrightarrow 0
$$
is the free resolution of $H$.  
Once one 
has the free resolution, one can read off the Betti 
numbers of $H$ as the multiplicities of the various free 
Hilbert modules that have appeared in the sequence, {\em in 
their order of appearance}.  
\end{rem}

We now discuss the examples of (1) for the case 
$r=3$ and arbitrary $d$.  

\begin{prop}
The 
Hilbert module associated with the 
$(d+3)$-contraction  $(S_1,\dots, S_d,0,0,0)$ 
acting on $H^2[z_1,\dots, z_d]$ 
has  Euler characteristic zero, and its sequence 
of Betti numbers is 
$$
(\beta_1,\dots, \beta_{d+3})=(1,3,3,1,0,\dots,0).
$$ 
\end{prop}

\begin{proof}[Sketch of Proof]
We show that the free resolution of $H$ 
has the form 
$$
0\xrightarrow{} F_4\xrightarrow{}
F_3
\xrightarrow{}
F_2
\xrightarrow{}
F_1\xrightarrow{} H\xrightarrow{} 0  
$$
where $F_k=H^2[z_1,\dots,z_{d+3}]\otimes G_k$, 
$G_1, G_2, G_3, G_4$ being graded coefficient spaces 
of respective dimensions $1, 3, 3, 1$.  We will 
exhibit the modules $F_k$ and the connecting maps 
explicitly, but we omit the details of computations 
with polynomials.  

We first compute the proper generator 
$H\ominus(Z\cdot H)$ of $H$.  Writing 
$$
T_1T_1^*+\cdots+T_{d+3}T_{d+3}^*=S_1S_1^*+\cdots+S_dS_d^*, 
$$
one sees that the defect operator 
$(\mathbf 1-\sum_k T_kT_k^*)^{1/2}$ is the one-dimensional 
projection $[1]$ onto the constant polynomials.  It follows 
that $H$ has defect $1$, and its proper 
generator is the one-dimensional space $\mathbb C\cdot 1$.  

Hence the first term in the free resolution of $H$ is given 
by the free cover $A_1: H^2[z_1,\dots,z_{d+3}]\to H$, where 
$A_1$ is the closure of the map defined on polynomials 
$f\in \mathbb C[z_1,\dots,z_{d+3}]$ by 
$$
A_1f=f(S_1,\dots,S_d,0,0,0)\cdot 1
=f(z_1,\dots,z_d,0,0,0).  
$$
$A_1$ is a coisometry, and further computation with 
polynomials shows that its kernel is 
the closure $K_1=\overline{(z_{d+1},z_{d+2},z_{d+3})}$ of 
the ideal in $\mathbb C[z_1,\dots,z_{d+3}]$ generated 
by $z_{d+1},z_{d+2},z_{d+3}$.  This gives a sequence 
of contractive homomorphisms  of degree zero 
$$
0\longrightarrow K_1
\longrightarrow 
H^2[z_1,\dots,z_{d+3}]\longrightarrow H
\longrightarrow 0  
$$  
that is exact at $H^2[z_1,\dots,z_{d+3}]$.

The kernel $K_1$ is a graded submodule of $H^2[z_1,\dots,z_{d+3}]$, but 
the rank of its defect operator is typically {\em infinite}.  
However, by Proposition \ref{propA}, 
it has a unique finite-dimensional 
proper generator $G$, given by 
$$
G=K_1\ominus (Z\cdot K_1)=K_1\ominus \overline{(z_1\cdot K_1+\cdots+z_{d+3}\cdot K_1)}.  
$$ 
To compute $G$, note that each of the elements 
$z_{d+1}$, $z_{d+2}$, $z_{d+3}$ is of degree one, while any 
homogeneous polynomial of $Z\cdot K_1$ is of degree at least two.  It 
follows that 
$K_1={\rm span}\{z_{d+1}, z_{d+2}, z_{d+3}\} \oplus (Z\cdot K_1)$,  
and this identifies $G$ as the $3$-dimensional 
Hilbert space 
$$
G={\rm span}\{ z_{d+1}, z_{d+2}, z_{d+3}\}.
$$
The multiplication 
map $A_2: F\otimes G\to F$ 
$$
A_2(f\otimes \zeta)=f\cdot \zeta,\qquad f\in\mathbb C[z_1,\dots,z_{d+3}], 
\quad \zeta\in G
$$
is a contractive morphism that defines a free cover 
of $K_1$; and $A_2$ becomes a 
degree zero map with respect to the gauge group $\Gamma$ 
on $H^2[z_1,\dots,z_{d+3}]\otimes G$ defined by 
$
\Gamma=\Gamma_0\otimes W
$
where $W$ is the restriction of the gauge group of 
$H^2[z_1,\dots,z_{d+3}]$ to its subspace $G$, namely 
$W(\lambda)=\lambda\mathbf 1_{G}$, $\lambda\in\mathbb T$.  
It follows that the sequence 
$$
H^2[z_1,\dots,z_{d+3}]\otimes G\underset{A_2}\longrightarrow 
H^2[z_1,\dots,z_{d+3}]
\underset{A_1}\longrightarrow H\longrightarrow 0
$$
is weakly exact at $H^2[z_1,\dots,z_{d+3}]$ and $H$.  

Now consider $K_2=\ker A_2\subseteq H^2[z_1,\dots,z_{d+3}]\otimes G$.  Since every 
element of $H^2[z_1,\dots,z_{d+3}]\otimes G$ can be written uniquely in the form 
$$
\xi_1\otimes z_{d+1}+\xi_2\otimes z_{d+2}+\xi_3\otimes z_{d+3}, 
\qquad \xi_k\in H^2[z_1,\dots,z_{d+3}]
$$
we have 
$$
K_2 = \{\xi_1\otimes z_{d+1}+\xi_2\otimes z_{d+2}+\xi_3\otimes z_{d+3}: 
z_{d+1}\cdot \xi_1+z_{d+2}\cdot \xi_2+z_{d+3}\cdot \xi_3=0\}.  
$$

A nontrivial calculation with polynomials 
now shows that $K_2$ is the 
closed submodule of 
of $H^2\otimes G$ generated by the three ``commutators" 
$\zeta_1, \zeta_2,\zeta_3$ 
\begin{align*}
\zeta_1&=z_{d+2}\otimes z_{d+3}-z_{d+3}\otimes z_{d+2}=z_{d+2}\wedge z_{d+3}, \\
\zeta_2&= z_{d+1}\otimes z_{d+3}-z_{d+3}\otimes z_{d+1}=z_{d+1}\wedge z_{d+3} \\
\zeta_3&=  z_{d+1}\otimes z_{d+2}-z_{d+2}\otimes z_{d+1}=z_{d+1}\wedge z_{d+2}.  
\end{align*}
Note, for example, that 
$$
f\cdot\zeta_1+g\cdot\zeta_2 = 
-gz_{d+3}\otimes z_{d+1} - fz_{d+3}\otimes z_{d+2}+(gz_{d+1}+fz_{d+2})\otimes z_{d+3}.
$$
These elements $\zeta_k$ 
are all homogeneous of degree two.  Since 
any homogeneous element of $Z\cdot K_2$ has degree at most 
three, it must be orthogonal to $\zeta_1, \zeta_2,\zeta_3$.  
It follows that 
$$
K_2\ominus (Z\cdot K_2)={\rm span}\{\zeta_2,\zeta_2,\zeta_3\}
$$
is $3$-dimensional, 
having $2^{-1/2}\zeta_1, 2^{-1/2}\zeta_2, 2^{-1/2}\zeta_3$
as an orthonormal basis.  

Set $\tilde G={\rm span}\{\zeta_2,\zeta_2,\zeta_3\}$, with its 
grading (in this case 
homogeneous of degree $2$) as inherited from the grading 
of $H^2[z_1,\dots,z_{d+3}]\otimes G$.  The corresponding free cover 
$A_3: H^2[z_1,\dots,z_{d+3}]\otimes \tilde G\to K_2$ is given by 
$$
A_3(f_1\otimes \zeta_1+f_2\otimes \zeta_2 + f_3\otimes \zeta_3)=
f_1\cdot\zeta_1+f_2\cdot\zeta_2+f_2\cdot\zeta_3, 
$$
for polynomials $f_1, f_2, f_3$, and the grading of 
$H^2[z_1,\dots,z_{d+3}]\otimes\tilde G$ is given by 
$\Gamma(\lambda)(f\otimes \zeta)=
\lambda^2(\Gamma_0(\lambda)f\otimes \zeta)$, 
$\lambda\in\mathbb T$.  

Finally, consider the submodule $K_3=\ker A_3\subseteq H^2[z_1,\dots,z_{d+3}]\otimes \tilde G$.  
Another computation with polynomials 
shows that $K_3$ has a single generator 
\begin{align*}
\eta &= z_{d-1}\otimes\zeta_1 - z_{d-2}\otimes  \zeta_2 + z_{d-3}\otimes\zeta_3 \\
&=
z_{d-1}\otimes (z_{d-2}\wedge z_{d-3}) - z_{d-2}\otimes (z_{d-1}\wedge z_{d-3}) +
z_{d-3}\otimes (z_{d-1}\wedge z_{d-2}), 
\end{align*}
where as above, $z_j\wedge z_k$ denotes 
$z_j\otimes z_k-z_k\otimes z_j$.  
The homogeneous element 
$\eta$ has degree $3$, so that after appropriate 
renormalization it becomes a unit vector spanning 
$K_3\ominus(Z\cdot K_3)$.  Thus, we obtain a free cover 
$A_4: H^2[z_1,\dots,z_{d+3}]\to K_3$ by closing the map of 
polynomials 
$$
A_4(f)=f\cdot\eta, \qquad f\in\mathbb C[z_1,\dots,z_{d+3}].  
$$ 
Notice that the grading that $H^2[z_1,\dots,z_{d+3}]$ acquires  
by this construction is not the standard grading $\Gamma_0$, 
but rather $\Gamma(\lambda)=\lambda^3\Gamma_0(\lambda)$, $\lambda\in \mathbb T$.  

Since the kernel of $A_4$ is obviously $\{0\}$, we have 
obtained a free resolution 
$$
0\xrightarrow{} F\xrightarrow{A_4}
F\otimes \tilde G
\xrightarrow{A_3}
F\otimes G
\xrightarrow{A_2}
F\xrightarrow{A_1} H\xrightarrow{} 0  
$$
in which $F=H^2[z_1,\dots,z_{d+3}]$.  

This shows
 that $H$ is a Hilbert module over $\mathbb C[z_1,\dots, z_{d+3}]$ 
whose Betti numbers $(\beta_1, \cdots, \beta_{d+3})$ are given 
by a nontrivial sequence 
$(1,3,3,1,0,\dots,0)$ with alternating sum zero.  
\end{proof}

\vfill

\bibliographystyle{alpha}

\newcommand{\noopsort}[1]{} \newcommand{\printfirst}[2]{#1}
  \newcommand{\singleletter}[1]{#1} \newcommand{\switchargs}[2]{#2#1}

\end{document}